\documentclass[12pt]{article}
\usepackage{amsmath,amsfonts,latexsym,amsthm,amssymb}
\newcommand{\F}{\mathbb F_q}
\newcommand{\K}{\overline{K}_c}
\newcommand{\R}{\mathcal R_K}
\newcommand{\RR}{\mathcal R_{\K}}
\newcommand{\A}{\mathcal A_K}
\newcommand{\AP}{\mathcal A_{K_{\text{perf}}}}
\newcommand{\KP}{K_{\text{perf}}}
\newcommand{\RP}{\mathcal R_{K_{\text{perf}}}}
\begin{document}
\newtheorem{prop}{Proposition}
\newtheorem{teo}{Theorem}
\pagestyle{plain}
\title{Strongly Nonlinear Differential Equations with Carlitz Derivatives 
over a Function Field}

\bigskip
\author{ANATOLY N. KOCHUBEI\footnote{This research was supported in part 
by CRDF under Grants UM1-2421-KV-02 and UM1-2567-OD-03.}\\ Institute of Mathematics,
National Academy\\ of Sciences of Ukraine, Tereshchenkivska 3, 
Kiev, 01601 Ukraine}
\date{}
\maketitle
\begin{abstract}
In earlier papers the author studied some classes of equations with Carlitz derivatives for 
$\F$-linear functions, which are the natural function field counterparts of linear ordinary
differential equations. Here we consider equations containing 
self-compositions $u\circ u\circ \cdots \circ u$ of the unknown 
function. As an algebraic background, imbeddings of the 
composition ring of $\F$-linear holomorphic functions into skew 
fields are considered.
\end{abstract}
\section{INTRODUCTION}

\medskip
Let $K$ be the set of formal Laurent series
$t=\sum\limits_{j=N}^\infty \xi_jx^j$ with coefficients $\xi_j$ 
from the Galois field $\F$, $\xi_N\ne 0$ if $t\ne 0$, $q=p^\upsilon $, 
$\upsilon \in \mathbf Z_+$,
where $p$ is a prime number. It is well known that $K$ is a 
locally compact field of characteristic $p$, with natural 
operations over power series, and the topology given by the 
absolute value $|t|=q^{-N}$, $|0|=0$. The element $x$ is a prime 
element of $K$. Any non-discrete locally compact field of characteristic 
$p$ is isomorphic to such $K$. Below we denote by $\K$ the completion 
of an algebraic closure $\overline{K}$ of $K$. The absolute value 
$|\cdot |$ can be extended in a unique way onto $\K$.

An important class of functions playing a significant part in the 
analysis over $\K$ is the class of $\F$-linear functions. A function 
$f$ defined on a $\F$-subspace $K_0$ of $K$ (or $\K$), with values
in $\K$, is called $\F$-{\it linear}
if $f(t_1+t_2)=f(t_1)+f(t_2)$ and $f(\alpha t)=\alpha f(t)$ for any 
$t,t_1,t_2\in K_0$, $\alpha \in \F$. A typical example is a 
$\F$-linear polynomial $\sum c_kt^{q^k}$ or, more generally, a 
power series $\sum\limits_{k=0}^\infty c_kt^{q^k}$, where $c_k\in 
\K$ and $|c_k|\le C^{q^k}$, convergent on a neighbourhood of the 
origin.

In the theory of differential equations over $K$ initiated in 
\cite{K2,K3} (which deals also with some non-analytic $\F$-linear 
functions) the role of a derivative is played by the operator
$$
d=\sqrt[q]{}\circ \Delta ,\quad (\Delta u)(t)=u(xt)-xu(t),
$$
introduced by Carlitz \cite{C1} and used subsequently in various 
problems of analysis in positive characteristic 
\cite{C2,G1,K1,Th,W}.

The differential equations considered so far were analogs of 
linear ordinary differential equations, though the operator $d$ is 
only $\F$-linear and the meaning of a polynomial coefficient in 
the function field case is not a usual multiplication by a 
polynomial, but the action of a polynomial in the $\F$-linear 
operator $\tau$, $\tau u=u^q$. Note that $\F$-linear polynomials 
form a ring with respect to the composition $u\circ v$ (the usual 
multiplication violates the $\F$-linearity), so that natural 
classes of equations with stronger nonlinearities must contain 
expressions like $u\circ u$ or, more generally, $u\circ u\cdots \circ 
u$. An investigation of such ``strongly nonlinear'' Carlitz 
differential equations is the main aim of this paper.

However we have to begin with algebraic preliminaries of some 
independent interest (so that not all the results are used in the 
subsequent sections) regarding the ring $\R$ of locally 
convergent $\F$-linear holomorphic functions. The ring is 
non-commutative, and the algebraic structures related to strongly 
nonlinear Carlitz differential equations are much more complicated 
than their classical counterparts. So far their understanding is 
only at its initial stage. Here we show that $\R$ is 
imbedded into a skew field of $\F$-linear ``meromorphic'' series 
containing terms like $t^{q^{-k}}$. Note that a deep investigation 
of bi-infinite series of this kind convergent on the whole of $\K$ 
has been carried out by Poonen \cite{Po}. We also prove an 
appropriate version of the implicit function theorem.

After the above preparations we consider general strongly 
nonlinear first order $\F$-linear differential equations (resolved 
with respect to the derivative of the unknown function) and prove 
an analog of the classical Cauchy theorem on the existence and 
uniqueness of a local holomorphic solution of the Cauchy problem. 
In our case the classical majorant approach (see e.g. \cite{Hi}) 
does not work, and the convergence is proved by direct estimates. 
We also consider a class of Riccati-type equations possessing 
$\F$-linear solutions which are meromorphic in the above sense.

\section{Skew fields of $\F$-linear power series}

\medskip
Let $\mathcal R_K$ be the set of all formal power series 
$a=\sum\limits_{k=0}^\infty a_kt^{q^k}$ where $a_k\in K$, 
$|a_k|\le A^{q^k}$, and $A$ is a positive constant depending on 
$a$. In fact each series $a=a(t)$ from $\R$ converges on a 
neighbourhood of the origin in $K$ (and $\K$).

$\R$ is a ring with respect to the termwise addition and the 
composition
$$
a\circ b=\sum\limits_{l=0}^\infty \left( \sum\limits_{n=0}^l
a_nb_{l-n}^{q^n}\right) t^{q^l},\quad b=\sum\limits_{k=0}^\infty 
b_kt^{q^k},
$$
as the operation of multiplication. Indeed, if $|b_k|\le B^{q^k}$, 
then, by the ultra-metric property of the absolute value,
$$
\left| \sum\limits_{n=0}^la_nb_{l-n}^{q^n}\right| \le 
\max\limits_{0\le n\le l}A^{q^n}\left( B^{q^{l-n}}\right)^{q^n}\le 
C^{q^l}
$$
where $C=B\max (A,1)$. The unit element in $\R$ is $a(t)=t$. It is 
easy to check that $\R$ has no zero divisors.

If $a\in \R$, $a=\sum\limits_{k=0}^\infty a_kt^{q^k}$, is such 
that $|a_0|\le 1$ and $|a_k|\le A^{q^k}$, $|A|\ge 1$, for all $k$, 
then we may write
$$
|a_k|\le A_1^{q^k-1},\quad k=0,1,2,\ldots ,
$$
if we take $A_1\ge A^{q^k/(q^k-1)}$ for all $k\ge 1$. If also 
$b=\sum\limits_{k=0}^\infty b_kt^{q^k}$, $|b_k|\le B_1^{q^k-1}$, 
$B_1\ge 1$, then for $a\circ b=\sum\limits_{l=0}^\infty 
c_lt^{q^l}$ we have
$$
|c_l|\le \max_{i+j=l}A_1^{q^i-1}\left( B_1^{q^j-1}\right)^{q^i}\le 
C_1^{q^l-1}
$$
where $C_1=\max (A_1,B_1)$. In particular, in this case the 
coefficients of the series for $a^n$ (the composition power) 
satisfy an estimate of this kind, with a constant independent of 
$n$.

\begin{prop}
The ring $\R$ is a left Ore ring, thus it possesses a classical 
ring of fractions.
\end{prop}

{\it Proof}. By Ore's theorem (see \cite{Her}) it suffices to show 
that for any elements $a,b\in \R$ there exist such elements 
$a',b'\in \R$ that $b'\ne 0$ and
\begin{equation}
a'\circ b=b'\circ a.
\end{equation}
We may assume that $a\ne 0$,
$$
a=\sum\limits_{k=m}^\infty a_kt^{q^k},\quad b=\sum\limits_{k=l}^\infty 
b_kt^{q^k},
$$
$m,l\ge 0$, $a_m\ne 0$, $b_l\ne 0$. 

Without restricting generality we may assume that $l=m$ (if we 
prove (1) for this case and if, for example, $l<m$, we set 
$b_1=t^{q^{m-l}}\circ b$, find $a'',b'$ in such a way that $a''\circ b_1=b'\circ 
a$, and then set $a'=a''\circ t^{q^{m-l}}$), and that 
$a_l=b_l=\alpha$, so that
$$
a=\alpha t^{q^l}+\sum\limits_{k=l+1}^\infty a_kt^{q^k},\quad 
b=\alpha t^{q^l}+\sum\limits_{k=l+1}^\infty b_kt^{q^k},
$$
$\alpha \ne 0$.

We seek $a',b'$ in the form
$$
a'=\sum\limits_{j=l}^\infty a_j't^{q^j},\quad b'=\sum\limits_{j=l}^\infty 
b_j't^{q^j}.
$$

The coefficients $a_j',b_j'$ can be defined inductively. Set 
$a_l'=b_l'=1$. If $a_j',b_j'$ have been determined for $l\le j\le 
k-1$, then $a_k',b_k'$ are determined from the equality of the 
$(k+l)$-th terms of the composition products:
$$
a_k'\alpha^{q^k}+\sum\limits_{\genfrac{}{}{0pt}{}{i+j=k+l}{j\ne 
l}}a_i'b_j^{q^i}=b_k'\alpha^{q^k}+\sum\limits_{\genfrac{}{}{0pt}{}{i+j=k+l}{j\ne 
l}}b_i'a_j^{q^i}
$$
(the above sums do not contain non-trivial terms with $a_i',b_i'$, 
$i\ge k$, since $a_j=b_j=0$ for $j<l$).

In particular, we may set $b_k'=0$,
$$
a_k'=\alpha^{-q^k}\left\{ \sum\limits_{\genfrac{}{}{0pt}{}{i+j=k+l}{i<k,j\ne 
l}}\left( a_i'b_j^{q^i}-b_i'a_j^{q^i}\right) \right\} .
$$
If this choice is made for each $k\ge l+1$, then we have $b_i'=0$ 
for every $i\ge l+1$, so that
\begin{equation}
a_k'=\alpha^{-q^k}\sum\limits_{\genfrac{}{}{0pt}{}{i+j=k+l}{i<k,j\ne 
l}}a_i'b_j^{q^i}.
\end{equation}

Denote $C_1=|\alpha |^{-1}$. We have $|b_j|\le C_2^{q^j}$ for all 
$j$. Denote, further, $C_3=\max (1,C_1,C_2)$, $C_4=C_3^{q^{l+2}}$. 
Let us prove that
$$
\left| a_k'\right| \le C_4^{q^k}.
$$

Suppose that $\left| a_i'\right| \le C_4^{q^i}$ for all $i$, $l\le i\le k-1$ 
(this is obvious for $i=1$, since $a_l'=1$). By (2),
\begin{multline*}
\left| a_k'\right| \le C_1^{q^k}\max \limits_{\genfrac{}{}{0pt}{}{i+j=k+l}{i<k,j\ne 
l}}C_4^{q^i}C_2^{q^{i+j}}\le C_1^{q^k}C_4^{q^{k-1}}C_2^{q^{k+l}}\\
\le C_3^{q^k+q^{k+l+1}+q^{k+l}}=C_3^{(1+q^l+q^{l+1})q^k}\le \left( 
C_3^{q^{l+2}}\right)^{q^k}=C_4^{q^k},
\end{multline*}
as desired. Thus $a'\in \R$. $\qquad \blacksquare$

\medskip
Every non-zero element of $\R$ is invertible in the ring of 
fractions $\A$, which is actually a skew field consisting of 
formal fractions $c^{-1}d$, $c,d\in \R$.

\begin{prop}
Each element $a=c^{-1}d\in \A$ can be represented in the form 
$a=t^{q^{-m}}a'$ where $t^{q^{-m}}$ is the inverse of $t^{q^m}$, 
$a'\in \R$.
\end{prop}

{\it Proof}. It is sufficient to prove that any non-zero element 
$c\in \R$ can be written as $c=c'\circ t^{q^m}$ where $c$ is 
invertible in $\R$.

Let $c=\sum\limits_{k=m}^\infty c_kt^{q^k}$, $c_m\ne 0$, $|c_k|\le 
C^{q^k}$. Then
$$
c=c_m\left( t+\sum\limits_{l=1}^\infty 
c_m^{-1}c_{m+l}t^{q^l}\right) \circ t^{q^m}
$$
where $\left| c_m^{-1}c_{m+l}\right| \le C_1^{q^l-1}$ for all 
$l\ge 1$, if $C_1$ is sufficiently large. Denote
$$
w=\sum\limits_{l=1}^\infty c_m^{-1}c_{m+l}t^{q^l},\quad 
c'=c_m(t+w).
$$

The series $(t+w)^{-1}=\sum\limits_{n=0}^\infty (-1)^nw^n$ 
converges in the standard non-Archimedean topology of formal power 
series (see \cite{Pierce}, Sect. 19.7) because the formal power 
series for $w^n$ begins from the term with $t^{q^m}$; recall that 
$w^n$ is the composition power, and $t$ is the unit element. 
Moreover, $w^n=\sum\limits_{j=n}^\infty a_j^{(n)}t^{q^j}$ where 
$\left| a_j^{(n)}\right| \le C_1^{q^j-1}$ for all j, with the same 
constant independent of $n$. Using the ultra-metric inequality we 
find that the coefficients of the formal power series 
$(t+w)^{-1}=\sum\limits_{j=0}^\infty a_jt^{q^j}$ (each of them is, 
up to a sign, a finite sum of the coefficients $a_j^{(n)}$) 
satisfy the same estimate. Therefore $(c')^{-1}\in \R$. $\qquad 
\blacksquare$

\medskip
The skew field of fractions $\A$ can be imbedded into wider skew 
fields where operations are more explicit. Let $\KP$ be the 
perfection of the field $K$. Denote by $\AP^\infty$ the 
composition ring of $\F$-linear formal Laurent series 
$a=\sum\limits_{k=m}^\infty a_kt^{q^k}$, $m\in \mathbb Z$, $a_k\in 
\KP$, $a_m\ne 0$ (if $a\ne 0$). Since $\tau$ is an automorphism of $\KP$,
$\AP^\infty$ is a special case of the well-known ring of twisted 
Laurent series \cite{Pierce}. Therefore $\AP^\infty$ is a skew 
field.

Let $\AP$ be a subring of $\AP^\infty$ consisting of formal series 
with $|a_k|\le A^{q^k}$ for all $k\ge 0$. Just as in the proof of 
Proposition 2, we show that $\AP$ is actually a skew field. Its 
elements can be written in the form $t^{q^{-m}}\circ c$ where $c$ 
is an invertible element of the ring $\RP \in \AP$ of formal power 
series $\sum\limits_{k=0}^\infty a_kt^{q^k}$. In contrast to the 
case of the skew field $\A$, in $\AP$ the multiplication of $t^{q^{-m}}$
by $c$ is indeed the composition of (locally defined) functions, 
so that $\AP$ consists of fractional power series understood in 
the classical sense.

Of course, $\AP$ can be extended further, by considering 
$\overline{K}$ or $\K$ instead of $\KP$. The above reasoning 
carries over to these cases (we can also consider the ring 
$\mathcal R_{\overline{K}_c}$ of locally convergent $\F$-linear 
power series as the initial ring). In each of them the presence of 
a fractional composition factor $t^{q^{-m}}$ is a $\F$-linear 
counterpart of a pole of the order $m$.

\section{Recurrent relations}

\medskip
In our investigations of strongly nonlinear equations and implicit 
functions we encounter recurrent relations of the same form
\begin{equation}
c_{i+1}=\mu_i\sum\limits_{\genfrac{}{}{0pt}{}{j+l=i}{l\ne 0}}
\sum\limits_{k=1}^\infty B_{jkl}\left( \sum\limits_{n_1+\cdots 
+n_k=l}c_{n_1}c_{n_2}^{q^{n_1}}\cdots c_{n_k}^{q^{n_1+\cdots 
+n_{k-1}}}\right)^{q^{j+\lambda}}+a_i,\quad i=1,2,\ldots ,
\end{equation}
(here and below $n_1,\ldots ,n_k\ge 1$ in the internal sum), with 
coefficients from $\K$, such that $|\mu_i|\le M$, $M>0$, $|B_{jkl}|\le B^{kq^j}$,
$B\ge 1$, $|a_i|\le M$ for all $i,j,k,l$; the number $\lambda$ is either 
equal to 1, or $\lambda =0$, and in that case $|B_{01l}|\le 1$.

\begin{prop}
For an arbitrary element $c_1\in \K$, the sequence determined by the 
relation (3) satisfies the estimate $|c_n|\le C^{q^n}$, $n=1,2,\ldots$, with
some constant $C\ge 1$.
\end{prop}

{\it Proof}. Set $c_n=\sigma d_n$, $|\sigma |<1$, $n=1,2,\ldots$, and 
substitute this into (3). We have
\begin{multline*}
d_{i+1}=\mu_i\Biggl\{ \sum\limits_{\genfrac{}{}{0pt}{}{j+l=i}{l\ne 0}}
\sum\limits_{k=1}^\infty B_{jkl}\sum\limits_{n_1+\cdots 
+n_k=l}
\sigma^{\left( 1+q^{n_1}+\cdots +q^{n_1+\cdots 
+n_{k-1}}\right)^{q^{j+\lambda }}-1}\\
\times \left( d_{n_1}d_{n_2}^{q^{n_1}}\cdots d_{n_k}^{q^{n_1+\cdots 
+n_{k-1}}}\right)^{q^{j+\lambda}}\Biggr\} +\sigma^{-1}a_i.
\end{multline*}
Here
$$
\left| \sigma^{\left(1+q^{n_1}+\cdots +q^{n_1+\cdots 
+n_{k-1}}\right)^{q^{j+\lambda }}-1}\right| \le |\sigma|^{kq^{j+\lambda }-1},
$$
and (under our assumptions) choosing such $\sigma$ that $|\sigma |$ is small enough
we reduce (3) to the relation
\begin{equation}
d_{i+1}=\mu_i\sum\limits_{\genfrac{}{}{0pt}{}{j+l=i}{l\ne 0}}
\sum\limits_{k=1}^\infty b_{jkl}\sum\limits_{n_1+\cdots 
+n_k=l}\left( d_{n_1}d_{n_2}^{q^{n_1}}\cdots d_{n_k}^{q^{n_1+\cdots 
+n_{k-1}}}\right)^{q^{j+\lambda}}+\sigma^{-1}a_i,\quad i=1,2,\ldots ,
\end{equation}
where $|b_{jkl}|\le 1$. 

It follows from (4) that
$$
|d_{i+1}|\le M\max\limits_{\genfrac{}{}{0pt}{}{j+l=i}{l\ne 0}}
\sup\limits_{k\ge 1}\max\limits_{n_1+\cdots +n_k=l}\max
\Biggl\{ \left( |d_{n_1}|\cdot |d_{n_2}^{q^{n_1}}|\cdots |d_{n_k}|^{q^{n_1+\cdots 
+n_{k-1}}}\right)^{q^{j+\lambda}},M^{-1}\left| \sigma^{-1}a_i\right| \Biggr\} .
$$

Let $B=\max \left\{ 
1,M,|d_1|,M^{-1}\sup\limits_i|\sigma^{-1}a_i|\right\}$. Let us 
show that
\begin{equation}
|d_n|\le B^{q^{n-1}+q^{n-2}+\cdots +1},\quad n=1,2,\ldots .
\end{equation}

This is obvious for $n=1$. Suppose that we have proved (5) for 
$n\le i$. Then
\begin{multline*}
|d_{i+1}|\le M\max\limits_{j+l=i}
\sup\limits_{k\ge 1}\max\limits_{n_1+\cdots +n_k=l}\left(
B^{q^{n_1-1}+q^{n_1-2}+\cdots +1}\cdot B^{q^{n_1+n_2-1}+q^{n_1+n_2-2}+\cdots 
+q^{n_1}}\cdots \right. \\ \left.
\times \cdots B^{q^{n_1+\cdots +n_{k-1}+n_k-1}+q^{n_1+\cdots 
+n_{k-1}+n_k-2}+\cdots +1}+q^{n_1+\cdots 
+n_{k-1}}\right)^{q^{j+1}}
\le M\max\limits_{j+l=i}B^{q^{j+l}+\cdots +q^{j+1}}\\
\le B\cdot B^{q^i+q^{i-1}+\cdots +q}=B^{q^i+q^{i-1}+\cdots +1},
\end{multline*}
and we have proved (5). Therefore
$$
|c_n|\le |\sigma |B^{\frac{q^n-1}{q-1}}\le C^{q^n}
$$
for some $C$, as desired. $\qquad \blacksquare$

\medskip
\section{Implicit functions of algebraic type}

\medskip
In this section we look for $\F$-linear locally holomorphic 
solutions of equations of the form
\begin{equation}
P_0(t)+P_1(t)\circ z+P_2(t)\circ (z\circ z)+\cdots +P_N(t)\circ 
\underbrace{(z\circ z\circ \cdots \circ z)}_N=0
\end{equation}
where $P_0,P_1,\ldots P_N\in \RR$. Suppose that the coefficient 
$P_k(t)=\sum\limits_{j\ge 0}a_{jk}t^{q^j}$ is such that 
$a_{00}=0$, $a_{01}\ne 0$; these assumptions are similar to the 
ones guaranteeing the existence and uniqueness of a solution in 
the classical complex analysis. Then (see Sect. 2) $P_1$ is 
invertible in $\RR$, and we can rewrite (6) in the form
\begin{equation}
z+Q_2(t)\circ (z\circ z)+\cdots +Q_N(t)\circ 
\underbrace{(z\circ z\circ \cdots \circ z)}_N=Q_0(t)
\end{equation}
where $Q_0,Q_2,\ldots ,Q_N\in \RR$, that is
$$
Q_k(t)=\sum\limits_{j=0}^\infty b_{jk}t^{q^j},\quad |b_{jk}|\le 
B_k^{q^j},
$$
for some constants $B_k>0$, and $b_{00}=0$.

\begin{prop}
The equation (6) has a unique solution $z\in \RR$ satisfying the 
``initial condition'' 
$$
\frac{z(t)}t\longrightarrow 0,\quad t\to 0.
$$
\end{prop}

\medskip
{\it Proof}. Let us look for a solution of the transformed 
equation (7), of the form
\begin{equation}
z(t)=\sum\limits_{i=1}^\infty c_it^{q^i},\quad c_i\in \K ;
\end{equation}
our initial condition is automatically satisfied for a function (8).

Substituting (8) into (7) we come to the system of equalities
\begin{equation} 
c_i=-\sum\limits_{k=2}^N\sum\limits_{\genfrac{}{}{0pt}{}{j+l=i}{j\ge 0,l\ge 1}}
b_{jk}\left( \sum\limits_{\genfrac{}{}{0pt}{}{n_1+\cdots +n_k=l}{n_j\ge 1}}
c_{n_1}c_{n_2}^{q^{n_1}}\cdots c_{n_k}^{q^{n_1+\cdots 
+n_{k-1}}}\right)^{q^j}+b_{i0},\quad i\ge 1.
\end{equation}
In each of them the right-hand side depends only on $c_1,\ldots 
,c_{i-1}$, so that the relations (9) determine the coefficients of 
a solution (8) uniquely. By Proposition 3, $z\in \RR$. $\qquad 
\blacksquare$

\medskip
More generally, let
$$
P_1(t)=\sum\limits_{j\ge \nu}a_{j1}t^{q^j},\quad \nu \ge 0,\quad 
a_{\nu 1}\ne 0.
$$
Then the equation (6) has a unique solution in $\RR$, of the form
$$
z(t)=\sum\limits_{i=\nu +1}^\infty c_it^{q^i},\quad c_i\in \K .
$$
The proof is similar.

\section{Equations with Carlitz derivatives}

\medskip
Let us consider the equation
\begin{equation}
dz(t)=\sum\limits_{j=0}^\infty \sum\limits_{k=1}^\infty 
a_{jk}\tau^j\underbrace{(z\circ z\circ \cdots \circ 
z)}_k(t)+\sum\limits_{j=0}^\infty a_{j0}t^{q^j}
\end{equation}
where $a_{jk}\in \K$, $|a_{jk}|\le A^{kq^j}$ ($k\ge 1$), 
$|a_{j0}|\le A^{q^j}$, $A\ge 1$. We look for a solution in the 
class of $\F$-linear locally holomorphic functions of the form
\begin{equation}
z(t)=\sum\limits_{k=1}^\infty c_kt^{q^k},\quad c_k\in \K ,
\end{equation}
thus assuming the initial condition $t^{-1}z(t)\to 0$, as $t\to 
0$.

\medskip
\begin{teo}
A solution (11) of the equation (10) exists with a non-zero radius 
of convergence, and is unique.
\end{teo}

\medskip
{\it Proof}. We may assume that
\begin{equation}
|a_{j0}|\le 1,\quad a_{j0}\to 0,\quad \text{as $j\to \infty$}.
\end{equation}
Indeed, if that is not satisfied, we can perform a time change 
$t=\gamma t_1$ obtaining an equation of the same form, but with 
the coefficients $a_{j0}\gamma^{q^j}$ instead of $a_{j0}$, and it 
remains to choose $\gamma$ with $|\gamma |$ small enough. Note 
that, in contrast with the case of the usual derivatives, the 
operator $d$ commutes with the above time change.

Assuming (12) we substitute (11) into (10) using the fact that 
$d\left( c_kt^{q^k}\right) =c_k^{1/q}[k]^{1/q}t^{q^{k-1}}$, $k\ge 
1$, where $[k]=x^{q^k}-x$. Comparing the coefficients we come to 
the recursion
$$
c_{i+1}=[i+1]^{-1}\sum\limits_{\genfrac{}{}{0pt}{}{j+l=i}{j\ge 0,l\ge 
1}}\sum\limits_{k=1}^\infty a_{jk}^q\left( \sum\limits_{n_1+\cdots +n_k=l}
c_{n_1}c_{n_2}^{q^{n_1}}\cdots c_{n_k}^{q^{n_1+\cdots 
+n_{k-1}}}\right)^{q^{j+1}}+a_{i0},\quad i\ge 1,
$$
where $c_1=[1]^{-1}a_{00}^q$. This already shows the uniqueness of 
a solution. The fact that $|c_i|\le C^{q^i}$ for some $C$ follows 
from Proposition 3. $\qquad \blacksquare$

\medskip
Using Proposition 4 we can easily reduce to the form (10) some 
classes of equations given in the form not resolved with respect 
to $dz$.

As in the classical case of equations over $\mathbb C$ (see 
\cite{Hi}), some of equations (10) can have also non-holomorphic 
solutions, in particular those which are meromorphic in the sense 
of Sect. 2. As an example, we consider Riccati-type equations
\begin{equation}
dy(t)=\lambda (y\circ y)(t)+(P(\tau )y)(t)+R(t)
\end{equation}
where $\lambda \in \K$, $0<|\lambda |\le q^{-1/q^2}$,
$$
(P(\tau )y)(t)=\sum\limits_{k=1}^\infty p_ky^{q^k}(t),\quad 
R(t)=\sum\limits_{k=0}^\infty r_kt^{q^k},
$$
$p_k,r_k\in \K$, $|p_k|\le q^{-1/q^2}$, $|r_k|\le q^{-1/q^2}$ for 
all $k$.

\medskip
\begin{teo}
Under the above assumptions, the equation (13) possesses solutions 
of the form
\begin{equation}
y(t)=ct^{1/q}+\sum\limits_{n=0}^\infty a_nt^{q^n},\quad c,a_n\in 
\K,\ c\ne 0,
\end{equation}
where the series converges on the open unit disk $|t|<1$.
\end{teo}

\medskip
{\it Proof}. For the function (14) we have
$$
dy(t)=c^{1/q}[-1]^{1/q}t^{q^{-2}}+\sum\limits_{n=1}^\infty 
a_n^{1/q}[n]^{1/q}t^{q^{n-1}},\quad [-1]=x^{1/q}-x,
$$
\begin{multline*}
(y\circ y)(t)=c\left( ct^{1/q}+\sum\limits_{n=0}^\infty 
a_nt^{q^n}\right)^{1/q}+\sum\limits_{n=0}^\infty a_n\left( 
ct^{1/q}+\sum\limits_{m=0}^\infty a_mt^{q^m}\right)^{q^n}\\
=c^{1+\frac{1}q}t^{q^{-2}}+\left( ca_0^{1/q}+ca_0\right) 
t^{q^{-1}}+\sum\limits_{n=0}^\infty \left( 
ca_{n+1}^{1/q}+c^{q^{n+1}}a_{n+1}\right) t^{q^n}
+\sum\limits_{l=0}^\infty t^{q^l}\sum\limits_{\genfrac{}{}{0pt}{}{m+n=l}{m,n\ge 
0}}a_na_m^{q^n}.
\end{multline*}
Finally,
$$
(P(\tau )y)(t)=\sum\limits_{k=0}^\infty p_{k+1}c^{q^{k+1}}t^{q^k}
+\sum\limits_{l=0}^\infty t^{q^l}\sum\limits_{\genfrac{}{}{0pt}{}{i+j=l}{i\ge 
1,j\ge 0}}p_ia_j^{q^i}.
$$

Comparing the coefficients we find that
\begin{equation}
c=\lambda^{-1}[-1]^{1/q},\quad a_0^{1/q}+a_0=0,
\end{equation}
\begin{equation}
a_{l+1}^{1/q}([l+1]^{1/q}-\lambda c)-\lambda 
c^{q^{l+1}}a_{l+1}=\lambda \sum\limits_{\genfrac{}{}{0pt}{}{m+n=l}{m,n\ge 
0}}a_na_m^{q^n}+\sum\limits_{\genfrac{}{}{0pt}{}{i+j=l}{i\ge 
1,j\ge 0}}p_ia_j^{q^i}+r_l,\quad l\ge 0.
\end{equation}

By (15), we have $|c|\ge 1$, and either $a_0=0$, or $|a_0|=1$. 
Next, (16) is a recurrence relation (with an algebraic equation to 
be solved at each step) giving values of $a_l$ for all $l\ge 1$. 
Let us prove that $|a_j|\le 1$ for all $j$. Suppose we have proved 
that for $j\le l$. It follows from (16) that
\begin{equation}
\left| 
a_{l+1}[l+1]-\lambda^qc^qa_{l+1}-\lambda^qc^{q^{l+2}}a_{l+1}^q\right| 
\le q^{-1/q}.
\end{equation}

Suppose that $|a_{l+1}|>1$. We have $\lambda^qc^q=[-1]$, so that 
$|\lambda^qc^q|=q^{-1/q}$, and since $|[l+1]|=q^{-1}$ and $|c|\ge 
1$, we find that 
$$
\left| a_{l+1}[l+1]\right| <\left| \lambda^qc^qa_{l+1}\right| 
<\left| \lambda^qc^{q^{l+2}}a_{l+1}^q\right| .
$$
Therefore the left-hand side of (17) equals $\left| \lambda^qc^q\right|
\cdot \left| c^{q^{l+1}}\right| \cdot \left| a_{l+1}^q\right| 
>q^{-1/q}$, and we have come to a contradiction. $\qquad 
\blacksquare$


\begin{thebibliography}{999}
\bibitem{C1}
L. Carlitz, On certain functions connected with polynomials in a
Galois field, {\it Duke Math. J.} {\bf 1} (1935), 137--168.
\bibitem{C2}
L. Carlitz, Some special functions over $GF(q,x)$, {\it Duke Math. J.} 
{\bf 27} (1960), 139--158.
\bibitem{G1}
D. Goss, Fourier series, measures, and divided power series in
the theory of function fields, {\it K-Theory} {\bf 1} (1989),
533--555.
\bibitem{Her}
I. N. Herstein, {\it Noncommutative Rings}, The Carus Math. 
Monograph No. 15, Math. Assoc. of America, J. Wiley and Sons, 
1968.
\bibitem{Hi}
E. Hille, {\it Lectures on Ordinary Differential Equations}, 
Addison-Wesley, Reading, 1969.
\bibitem{K1}
A. N. Kochubei, $\F$-linear calculus over function fields, {\it J.
Number Theory} {\bf 76} (1999), 281--300.
\bibitem{K2}
A. N. Kochubei, Differential equations for $\F$-linear functions, {\it J.
Number Theory} {\bf 83} (2000), 137--154.
\bibitem{K3}
A. N. Kochubei, Differential equations for $\F$-linear functions,
II: Regular singularity, {\it Finite Fields Appl.} {\bf 9} (2003), 
250--266.
\bibitem{Pierce}
R. S. Pierce, {\it Associative Algebras}, Springer, New York, 
1982.
\bibitem{Po}
B. Poonen, Fractional power series and pairings on Drinfeld 
modules, {\it J. Amer. Math. Soc.} {\bf 9} (1996), 783--812.
\bibitem{Th}
D. Thakur, Hypergeometric functions for function fields II,
{\it J. Ramanujan Math. Soc.} {\bf 15} (2000), 43--52.
\bibitem{W}
C. G. Wagner, Linear operators in local fields of prime
characteristic, {\it J. Reine Angew. Math.} {\bf 251} (1971), 153--160.

\end{thebibliography}
\end{document}